\documentclass{article}

\usepackage{mathtools}
\usepackage{amsmath}
\usepackage{amsfonts}
\usepackage{amsthm}
\newtheorem{definition}{Definition}

\author{Jacob C. Vandenberg\thanks{School of Mathematical Sciences, Monash University, Clayton, Victoria 3800, Australia.}
\and Mark B. Flegg\thanks{School of Mathematical Sciences, Monash University, Clayton, Victoria 3800, Australia.}}

\title{Turing pattern or system heterogeneity? A numerical continuation approach to assessing the role of Turing instabilities in heterogeneous reaction-diffusion systems}

\newcommand{\basex}[1]{\mathbf{u}^{\star}_{#1}(\mathbf{x})}
\newcommand{\base}[1]{\mathbf{u}^{\star}_{#1}}

\DeclareOldFontCommand{\rm}{\normalfont\rmfamily}{\mathrm}
\DeclareOldFontCommand{\sf}{\normalfont\sffamily}{\mathsf}
\DeclareOldFontCommand{\tt}{\normalfont\ttfamily}{\mathtt}
\DeclareOldFontCommand{\bf}{\normalfont\bfseries}{\mathbf}
\DeclareOldFontCommand{\it}{\normalfont\itshape}{\mathit}

\newcommand{\OO}{\mathcal{O}}

\begin{document}

\maketitle

\begin{abstract}
Turing patterns in reaction-diffusion (RD) systems have classically been studied only in RD systems which do not explicitly depend on independent variables such as space. In practise, many systems for which Turing patterning is important are not homogeneous with ideal boundary conditions. In heterogeneous systems with stable steady states, the steady states are also necessarily heterogeneous which is problematic for applying the classical analysis. Whilst there has been some work done to extend Turing analysis to some heterogeneous systems, for many systems it is still difficult to determine if a stable patterned state is driven purely by system heterogeneity or if a Turing instability is playing a role. In this work, we try to define a framework which uses numerical continuation to map heterogeneous RD systems onto a sensible nearby homogeneous system. This framework may be used for discussing the role of Turing instabilities in establishing patterns in heterogeneous RD systems. We study the Schnakenberg and Gierer-Meinhardt models with spatially heterogeneous production as test problems. It is shown that for sufficiently large system heterogeneity (large amplitude spatial variations in morphogen production) it is possible that Turing-patterned and base states become coincident and therefore impossible to distinguish. Other exotic behaviour is also shown to be possible. We also study a novel scenario in which morphogen is produced locally at levels that could support Turing patterning but on intervals/patches which are on the scale of classical critical domain lengths. Without classical domain boundaries, Turing patterns are allowed to bleed through; an effect noted by other authors. In this case, this phenomena effectively changes the critical domain length. Indeed, we even note that this phenomena may also effectively couple local patches together and drive instability in this way. 
\end{abstract}



\section{Introduction}

    The reaction-diffusion (RD) equation is a nonlinear partial differential equation which exhibits extraordinary diverse behavior observed particularly in the life sciences \cite{Timm1992, SUN201643, doi:10.1126/science.269.5222.331}.
    It models the concentration of different species in time as they interact whilst diffusing in space relative to each other.
    The \emph{species} of the system could refer to a chemical species, biological species or ecological species, amongst other possibilities~\cite{Mndez2010ReactionTransportSM}.

    Under certain conditions, solutions to the RD equation can have an instability which is ``driven by diffusion''.
    This is called a Turing instability, which is usually defined as follows.
    \textit{Turing instabilities occur when an RD system has a spatially-uniform steady state which is unstable in the presence of diffusion, but stable in the absence of diffusion.}
    Alan Turing's seminal paper analyses Turing instabilities as a mechanism for explaining the emergence of spatial heterogeneity in diffuse biological chemical systems~\cite{Turing1952}.
    The reason Turing instabilities can explain this onset of heterogeneity is because they typically produce \emph{Turing patterns}.
    Turing patterns are stable solutions to the RD equation which have large spatial oscillations, and are stationary in time.
    Usually diffusion has the effect of ``flattening'' the solution.
    In this case, however, diffusion is what causes the system to deviate away from uniformity.
    
    Often, RD models are spatially homogeneous in the sense that the RD PDE does not explicitly contain the spatial variable $x$ (or $t$). Typically, RD models which exhibit Turing patterning are studied as homogeneous systems to simplify the analysis of the PDE (finding steady states, performing linear stability analysis, demonstrating the potential for patterning etc.). At the same time, most real world applications almost certainly contain spatial variation in model parameters.  
    Consider, for example, the patterning and development of digits, kidneys and lungs where homogeneous models are analysed for the presence of Turing instabilities despite there being obvious spatial heterogeneity in morphogen production rates ~\cite{LAWSON201643,doi:10.1126/science.1226804}.

    Understanding Turing patterning in the presence of spatially heterogeneous RD PDEs is not well understood and surprisingly has received very little attention in the literature. Perhaps, one of the reasons for this is that Turing analysis of spatially heterogeneous RD PDEs is challenging as it is not even necessarily apparent even how Turing instabilities should be defined. To begin, the unstable uniform steady state required for defining the Turing instability does not exist by definition for spatially heterogeneous RD PDEs. 
    
    The analysis by Krause et al. presents a general stability theory for a heterogeneous RD PDE. This paper is however limited to cases where heterogeneity varies slowly almost everywhere relative to the domain size~\cite{Krause2020}.
    In the paper, Krause et al. define a `base state' solution which replaces the notion of the uniform steady state which has been `flattened' by diffusion. The base state, which must be a stationary solution to the PDE, has certain properties. 
    Importantly, the base state does not have spatial oscillations with periods much smaller than the inhomogeneity in the PDE (it is nice and `diffused').
    Aside from this definition being vague, it is not clear that it should be the case if the PDE contains heterogeneities which vary on the same spatial scale as the Turing patterns for the system. This is because it is not easy to distinguish between patterned and base states if oscillations in the patterned state are on the same spatial scale as the base state. We shall also be adopting the term `base state' but attempting to find a more general approach to finding it.

    Another method which has been widely used in the literature is to limit the scope of the study to more specific examples.
    This includes choosing specific reaction terms such that an exact solution can be computed~\cite{PAGE200380, Auchmuty1975}.
    At this point, a stability analysis similar to the classical analysis can be performed.
    Using a linear reaction term is common~\cite{PAGE200380, PhysRevE.100.042220}, but nonlinear reaction terms can also be considered~\cite{Auchmuty1975}.
    Truncated Galerkin expansions of the solution have been used to study the stability of heterogeneous problems~\cite{VanGorder2021, PhysRevE.100.042220}.
    These too use specific examples to find base states analytically.
    No insight is given as to why the solutions that were found should be analogous to the uniform base state in the homogeneous case.

    In this manuscript, our aim is to investigate a method which may be used to find base states for heterogeneous reaction-diffusion PDEs. The stability of these base states may be used to define Turing patterns. We propose a method for describing base states and apply this method to the canonical Schnakenberg (substrate depletion) system as well as the Gierer-Meinhardt (activator-inhibitor) system. In both of these systems we allow the production of species to vary in space. We focus on two main curiosities. The first deals with critical phenomena which place limitations on when a base state may be defined and the second deals with the onset of critical domain lengths for Turing instabilities in the presence of heterogeneous production.

\section{Methods}

The classical spatially-homogeneous dimensionless reaction-diffusion system is
\begin{align}
    \frac{\partial\mathbf{u}}{\partial t} &= \mathsf{D}\nabla^2\mathbf{u} + \gamma\mathbf{F}(\mathbf{u}),~ \text{on}~\Omega, \label{Methods:HomRD}\\
    \nabla \mathbf{u} \cdot \mathbf{n} &= \mathbf{0},~ \text{on}~\partial\Omega.\label{Methods:HomBC}
\end{align}
Here $\mathbf{u}$ is a vector containing the concentration of model species/chemicals, $\mathsf{D}$ is a diagonal matrix of diffusion constants (with $\mathsf{D}_{11} = 1$ providing a characteristic timescale for nondimensionalisation), and $\mathbf{F}$ is a nonlinear vector-valued function describing the possible sources and sinks of, and reactions between, the species. The domain $\Omega$ (which has an outward normal vector $\mathbf{n}$) has been scaled through non-dimensionalisation so that the spatial scale of the system relative to that of diffusion is described by the magnitude of $\gamma$.

A Turing analysis of this system begins by finding the uniform steady state solution \(\mathbf{u}^{\star}\) such that \(\mathbf{F}(\mathbf{u}^{\star}) = \mathbf{0}\). Indeed this uniform state is a solution to the model because derivatives of $\mathbf{u}^{\star}$ (a constant) is zero.
Subsequently, a Turing pattern is formed when the solution \(\mathbf{u}^{\star}\), which is stable if $\mathsf{D} = \mathsf{0}$, is unstable. The uniform solution to the model \(\mathbf{u}^{\star}\) will be called the base state and in heterogeneous problems loses its uniformity. This is the natural, diffusion-flattened, state of the system.

We can extend the RD model to account for explicit spatial variation
\begin{align}
    \frac{\partial \mathbf{u}}{\partial t} &= \text{div}(\mathsf{D}(\mathbf{x})\nabla\mathbf{u}) + \gamma\mathbf{F}(\mathbf{u},\mathbf{x}),~ \text{on}~\Omega,\label{Ch1:eq:HetRD}\\
    \nabla\mathbf{u} \cdot \mathbf{n} &= \mathbf{0}~\text{on}~\partial \Omega.\label{Ch1:eq:HetBC}
\end{align}
If we were to proceed as before, we can take \(\mathbf{u}^{\star}(\mathbf{x})\) which satisfies \(\mathbf{F}(\mathbf{u}^{\star}(\mathbf{x}), \mathbf{x}) = \boldsymbol{0}\) for all \(\mathbf{x} \in \Omega\).
The diffusion term \(\text{div}(\mathsf{D}(\mathbf{x})\nabla\mathbf{u}^{\star}(\mathbf{x}))\) is not zero in general, which would mean \(\mathbf{u}^{\star}(\mathbf{x})\) is not a steady state solution of Equation (\ref{Ch1:eq:HetRD}).
Thus, it does not make sense to analyse its stability.
So in order to extend the definition of a Turing instability, we need to find a different base state \(\mathbf{u}^{\star}(\mathbf{x})\) which satisfies the steady state problem for Equations (\ref{Ch1:eq:HetRD}) and (\ref{Ch1:eq:HetBC}) but also should not be called a Turing pattern. Whilst a `pattern' is often defined as any stable stationary heterogeneous solution, we reserve the definition of pattern in this manuscript to describe any stationary heterogeneous state separate to the base state.

As it stands, there is no conventional way of finding or defining more generally what this base state is.
The only thing that can be said about the base state \(\mathbf{u}^{\star}(\mathbf{x})\) is that it should be somehow sensibly analogous to the uniform base state described for the homogeneous system. 

 We will narrow the scope of our efforts to investigate this system to the case where heterogeneity is in the reaction term only. Specifically, we look at systems with heterogeneous production rates of each species as we believe that this system is ubiquitous in biological application where morphogen is deferentially expressed in space but reactions between morphogens are autonomous as one might expect.
Thus, the form of the RD equation that we will be analysing is as follows and splits $\mathbf{F}$ up into autonomous, homogeneous $\mathbf{\hat{F}}$ and heterogeneous $\mathbf{G}$ components. How this partition should be done appropriately and uniquely we will discuss here, outlining the approach that we have taken, but we will justify this approach in Section \ref{BaseStates}.
\begin{align}
    \frac{\partial \mathbf{u}}{\partial t} &= \mathsf{D}\nabla^2\mathbf{u} + \gamma\left(\mathbf{\hat{F}}(\mathbf{u}) + \mathbf{G}(\mathbf{u},\mathbf{x}) \right),~ \text{on}~\Omega,\label{Ch1:parameterisedfull}\\
    \nabla\mathbf{u} \cdot \mathbf{n} &= \mathbf{0}~\text{on}~\partial \Omega.\label{Ch1:parameterisedBCfull}
\end{align}
To analyse this system, we will find it useful to `grow' the heterogeneous components by means of a parameter $\theta$ by defining the parameterised problem
\begin{align}
    \frac{\partial \mathbf{u}}{\partial t} &= \mathsf{D}\nabla^2\mathbf{u} + \gamma\left(\mathbf{\hat{F}}(\mathbf{u}) + \theta \mathbf{G}(\mathbf{u},\mathbf{x}) \right),~ \text{on}~\Omega,\label{Ch1:parameterised}\\
    \label{Ch1:parameterisedBC} \nabla\mathbf{u} \cdot \mathbf{n} &= \mathbf{0}~\text{on}~\partial \Omega.
\end{align}
Importantly, the parameter $\theta$ in these models describe the amplitude of the heterogeneity in the system and when $\theta \rightarrow 0$ a classical system is recovered and when $\theta \rightarrow 1$ the full heterogeneous problem is recovered. Importantly, as $\theta$ may be thought of as the amplitude of the heterogeneity and easily absorbed into $\mathbf{G}$, it is possible to also think of $\theta$ growing beyond 1 and simply forming part of a growing $\mathbf{G}$ in Equations (\ref{Ch1:parameterisedfull}) and (\ref{Ch1:parameterisedBCfull}).

Whilst there is freedom in the choice of the partition of $\mathbf{F}$ in Equation (\ref{Ch1:eq:HetRD}) into $\mathbf{G}$ and $\mathbf{\hat{F}}$ in Equation (\ref{Ch1:parameterisedfull}), we find it appropriate to uniquely define  $\mathbf{G}$ and $\mathbf{\hat{F}}$ for a given $\mathbf{F}$ in the following way.
\begin{align}
    \mathbf{\hat{F}} &= \frac{1}{|\Omega|} \int_{\Omega}\mathbf{F}(\mathbf{u},\mathbf{x})\,\partial \mathbf{x},\label{Ch1:choiceF}\\
    \mathbf{G} &= \mathbf{F} - \mathbf{\hat{F}}.\label{Ch1:choiceG}
\end{align}
This is a convenient choice when the reaction term can be decomposed into a spatially-independent coupling term and a spatially-dependent source term, resulting in the following.
\begin{align*}
    \mathbf{F}(\mathbf{u}, \mathbf{x}) = \mathbf{\hat{F}}(\mathbf{u}) + \mathbf{G}(\mathbf{x}),
\end{align*}
where the average value of \(\mathbf{G}\) is 0. Furthermore, by using this decomposition for $\mathbf{F}$, we ensure that for each $\theta$ the parameterised system (Equation (\ref{Ch1:parameterised})) adheres to the same decomposition rules whilst at the same time capturing autonomous reactions in $\mathbf{F}$ within $\mathbf{\hat{F}}$ and often it is these terms which are the characteristically important ingredients in the Turing behaviour of the system (noting that $\mathbf{F}\rightarrow \mathbf{\hat{F}}$ as $\theta \rightarrow 0$).

\subsection{Base states}\label{BaseStates}

In this section, we attempt to redefine the base state of a heterogeneous reaction-diffusion system as a parameterised continuation of a nearby homogeneous system. A necessary condition on the base state of a reaction-diffusion system (Equations (\ref{Ch1:parameterised}) and (\ref{Ch1:parameterisedBC})) is that it must be a stationary solution, against which stability can be later checked.

The base state of Equations (\ref{Ch1:parameterised}) and (\ref{Ch1:parameterisedBC}) shall be labelled as \(\basex{\theta}\) (and sometimes as \(\base{\theta}(\mathbf{x};\theta)\) to highlight dependence on the parameter $\theta$). We have that \(\basex{\theta}\) is a solution to 
\begin{align}
    \mathsf{D}\nabla^2\mathbf{u} + \gamma\left(\mathbf{\hat{F}}(\mathbf{u}) + \theta \mathbf{G}(\mathbf{u},\mathbf{x}) \right) &= \mathbf{0}~\text{on}~ \Omega, \label{MEQ1}\\
    \nabla\mathbf{u} \cdot \mathbf{n} &= \mathbf{0}~\text{on}~\partial \Omega. \label{MEQ1BC}
\end{align}
Since the base state should become the uniform steady state as $\theta\rightarrow 0$, we have that \(\base{0} \in \mathbb{R}^{N_s}\) (where $N_s$ is the number of species in the model) is constant in $\mathbf{x}$ and $\mathbf{\hat{F}}(\base{0}) = \mathbf{0}$. 

It makes sense to represent Equations (\ref{MEQ1}) and (\ref{MEQ1BC}) as the single equation
\begin{align}\label{baseCombined}
    \Phi(\mathbf{u},\mathbf{x},\mathbf{\bar{x}};\theta) = \left( \mathsf{D}\nabla^2\mathbf{u}(\mathbf{x}) + \gamma\left(\mathbf{\hat{F}}(\mathbf{u}(\mathbf{x})) + \theta \mathbf{G}(\mathbf{u}(\mathbf{x}),\mathbf{x})\right), \nabla\mathbf{u}(\mathbf{\bar{x}}) \cdot \mathbf{n}\right)^{\top}=\mathbf{0},
\end{align}
where $\mathbf{x}\in\Omega$ and $\mathbf{\bar{x}}\in\partial \Omega$.

In order to label a solution to \ref{baseCombined} as a base solution, we will further require that it varies continuously with respect to $\theta$. In this way, the base states of the system are tied, via continuation of the parameter $\theta$, to the base state $\base{0}$ (uniform steady state) of an associated homogeneous system (as $\theta\rightarrow 0$). 

\iftrue
    To ensure the existence of $\base{\theta}$ for some \(\theta \ne 0\), we can find some \(\eta > 0\) and \(\base{\theta} \::\:(-\eta,\eta)\to C^2(\Omega, \mathbb{R}^{N_s})\) such that \(\base{\theta}\) uniquely solves~\ref{baseCombined} and $\base{0}$ solves $\mathbf{\hat{F}}(\base{0}) = \mathbf{0}$. The value \(\eta\) provides a region where any \(-\eta \le \theta \le \eta\) is guaranteed to have a base state solution. Outside of \((-\eta,\eta)\), the amplitude of the heterogeneity may become so large that it is not possible to draw a continuation from $\base{0}$.
\fi

We define the Jacobian 
\begin{align}
    \mathbf{\bar{J}}_\theta(\mathbf{u},\mathbf{x}) &= \frac{\partial \Phi}{\partial \mathbf{u}} = \left( \mathbf{J}_\theta(\mathbf{u},\mathbf{x}), \mathbf{n} \cdot \nabla\right)^{\top} \\ &= 
    \left( \mathsf{D}\nabla^2 + \gamma \left(\mathbf{j}_{\hat{\mathbf{F}}}(\mathbf{u}) + \theta\mathbf{j}_{\mathbf{G}}(\mathbf{u},\mathbf{x}) \right), \mathbf{n} \cdot \nabla\right)^{\top}. \label{Jacobiandef}
\end{align}

Here $\mathbf{j}_{\hat{\mathbf{F}}}(\mathbf{u})$ and  $\mathbf{j}_{\mathbf{G}}(\mathbf{u},\mathbf{x}) $ are the Jacobians of $\hat{\mathbf{F}}$ and $\mathbf{G}$ respectively. For continuity and uniqueness of $\base{\theta}$ in $\theta$ 
at $\theta=0$ by the Implicit Function Theorem (IFT)~\cite{Chow2011}, we require $ \mathbf{\bar{J}}_0$ to reversible at $\base{0}$ and therefore, we require that $\mathbf{j}_{\hat{\mathbf{F}}}(\base{0})$ is nonsingular.




 Singularity in $\bar{\mathbf{J}}_\theta$ allows for the possibility that $\theta$ may become too large in magnitude for there to exist a defined base state $\base{\theta}$. It's unclear in general how large a heterogeneity ($\theta$) can get before the base state either stops existing or is not unique, or even if the base state is even bound in this way at all. Defining the base state outside of some potential maximum range $\theta^-<\theta<\theta^+$, is problematic and in our framework not (yet) possible. The values of $\theta^-$ and $\theta^+$ coincide with folds in the solution to Equations (\ref{MEQ1}) and (\ref{MEQ1BC}) characterised by singularities in $\bar{\mathbf{J}}_{\theta^-}$ and $\bar{\mathbf{J}}_{\theta^+}$.

\begin{definition}[Spatially-dependent Turing base state]
    \label{hetBaseState}
    For each \(\mathbf{u}_0 \in \mathbb{R}^{N_s}\) such that \(\mathbf{\hat{F}}(\mathbf{u}_0) = \mathbf{0}\) we define the associated spatially-dependent Turing base state (or just base state) for Equation (\ref{Ch1:parameterisedfull}) as follows.
    If there exists \(\base{\theta}(\mathbf{x};\theta) \in C^1(\Omega \times (0,1],\mathbb{R}^{N_s})\) which is a steady state solution to Equation (\ref{Ch1:parameterised}) for all \(\theta \in (0,1]\) and where  \(\base{0}(\mathbf{x};0) = \mathbf{u}_{0}\).
    Then \(\basex{1}\) is a Turing base state  to the spatially-dependent RD system~(Equations (\ref{Ch1:parameterisedfull}) and (\ref{Ch1:parameterisedBCfull})) associated with the uniform base state \(\basex{0}\).
\end{definition}

Defining the base state in this way is a natural extension of the classical homogeneous case, since the heterogeneous base state should not deviate too far from the uniform one in the situation where the amplitude of the heterogeneity in the system is small. In other words, if heterogeneity in the system is small, we would expect that the base state should be \textit{almost} `flat' from diffusion.

As an important note, we have chosen to define $\mathbf{\hat{F}}$ and $\mathbf{G}$ using Equations (\ref{Ch1:choiceF}) and (\ref{Ch1:choiceG}), in doing so we ensure that all autonomous terms in $\mathbf{F}$ (for example reaction kinetics between species which drive Turing instabilities) are encapsulated in $\mathbf{\hat{F}}$. Clearly, it is possible to simply define $\mathbf{G} = \mathbf{F}$ and $\mathbf{\hat{F}}=0$. With this choice, we immediately see that $\mathbf{j}_{\hat{\mathbf{F}}}(\base{0})$ is singular and continuation to the heterogeneous base state is impossible.

In the case where \(\mathbf{\hat{F}} \ne \mathbf{0}\),
we have
\begin{align*}
    \mathbf{\bar{J}}_0(\base{0},x)= \left.\frac{\partial \Phi}{\partial \mathbf{u}}\right|_{\mathbf{u} = \base{0},\theta = 0} = \left( \mathsf{D}\nabla^2 + \gamma\mathsf{J}_{\mathbf{\hat{F}}}(\base{0}), \mathbf{n} \cdot \nabla\right)^{\top}.
\end{align*}
We apply this to \(\mathbf{c}_j \hat{w}_m\), where \(\mathbf{c}_j \in \mathbb{R}^{N_s}\) is the \(j\)th eigenvector of \(\mathsf{A}_m = -\mathsf{D}k_m^2 + \gamma\mathsf{J}_{\mathbf{\hat{F}}}(\base{0})\) and $\hat{w}_m$ is the eigenfunction solving $\nabla^2 \hat{w}_m = -k_m^2\hat{w}_m $ on $\Omega$ with $ \nabla \hat{w}_m \cdot  \mathbf{n}= 0$ on $\partial \Omega$.
This gives us the following,
\begin{align*}
    \left( \mathsf{D}\nabla^2\mathbf{c}_j \hat{w}_m + \gamma\mathsf{J}_{\mathbf{\hat{F}}}\mathbf{c}_j \hat{w}_m, \mathbf{n} \cdot \nabla\mathbf{c}_j \hat{w}_m\right)^{\top} &= \left(\mathsf{A}_m \mathbf{c}_j \hat{w}_m, 0\right)^{\top}\\
    & = \left(\lambda_{j}(\mathsf{A}_m), 0\right)^{\top}\mathbf{c}_j\hat{w}_m, 
\end{align*}
where \(\lambda_{j}(\mathsf{A}_m)\) is the eigenvalue associated with the eigenvector \(\mathbf{c}_j\).
This eigenvalue determines the stability of the eigenvector \(\mathbf{c}_j \hat{w}_m\).
So if any eigenvector \(\mathbf{c}_j\) has a corresponding \(\lambda_{j}(\mathsf{A}_m) = 0\), the operator \(\mathbf{\bar{J}}_0 \) will not be invertible and the conditions for the IFT would not be satisfied.

The continuation of base states from $\theta = 0$ cannot proceed unless \(\mathbf{G}(\base{0}, \mathbf{x})\) is orthogonal to every eigenvector in the null space of the adjoint operator \(\frac{\partial \Phi}{\partial \mathbf{u}}^{*}\). That is,
\begin{align*}
    \int_{\Omega} \mathbf{G}(\base{0}, \mathbf{x})^{\top} \mathbf{v}\,\mathrm{d}\mathbf{x} = 0,~\forall \mathbf{v}\in \text{null} \left( \mathsf{D}\nabla^2+ \gamma\mathsf{J}^{\top}_{\mathbf{\hat{F}}}\right).
\end{align*}
This is a result of Fredholm's alternative~\cite{Brezis2011}. This solvability condition is not guaranteed.
So for any chosen parameterisation, there may still be cases where continuation is impossible about \(\theta = 0\).

We have chosen to multiply the heterogeneity \(\mathbf{G}\) by a parameter \(\theta\). Of course, this parameterisation of heterogeneity ($\theta = 1$) from the associated homogeneous system ($\theta = 0$) is not unique.
In Equations (\ref{Ch1:parameterised}) and (\ref{Ch1:parameterisedBC}) we increase the size of the heterogeneity linearly with the parameter $\theta$.
A more general parameterisation could be 
\begin{align}
    \frac{\partial \mathbf{u}}{\partial t} = \mathsf{D}\nabla^2\mathbf{u} + \gamma\mathbf{\hat{F}}(\mathbf{u}) + \gamma\mathbf{G}(\mathbf{u},\mathbf{x}; \theta),\label{Ch1:parameterised2}
\end{align}
provided that \(\mathbf{\hat{F}}(\mathbf{u}) + \mathbf{G}(\mathbf{u},\mathbf{x};1) \equiv \mathbf{F}(\mathbf{u},\mathbf{x})\), and \(\mathbf{G}(\mathbf{u},\mathbf{x}; 0) \equiv \mathbf{0}\).

The IFT only provides information about the existence and uniqueness of the base state solution branch locally.
The existence and uniqueness of the base state solution at \(\theta = 1\) is unknown \textit{a priori}.
In particular, it is unknown whether changing the parameterisation of \(\mathbf{G}\) will lead to a change in the base state or the existence of the base state.
For this, a global homotopy result would be required.

The analysis by Krause et al. gives general stability theory for a large perturbation in the limit as \(\gamma\) approaches \(\infty\)~\cite{Krause2020}.
However, little attention is given on redefining the base state for the Turing instability.
The analysis assumes that a steady state solution to the full RD equation~(Equations (\ref{Ch1:eq:HetRD}) and (\ref{Ch1:eq:HetBC})) exists, and that this solution has certain properties.
The first property is that the solution does not have spatial oscillations on the scale \(\OO(1/\epsilon)\).
This is an \textit{a posteriori} assumption, since no method is provided for determining whether the base state \(\mathbf{u}^{\star}(\mathbf{x})\) has \(\OO(1/\epsilon)\) oscillations without first finding \(\mathbf{u}^{\star}(\mathbf{x})\).
Since the heterogeneous RD equation is nonlinear in general, finding such a solution is non-trivial.
Finally, it is assumed that \(\mathbf{F}\) satisfies the boundary conditions \(\frac{\partial\mathbf{u}}{\partial x} = \mathbf{0}\) at \(x = 0,1\).

\subsection{Case studies}
In our numerical investigation, we focus attention on two popular models; the Schnakenberg model and the Gierer-Meinhardt model. In their standard homogeneous forms, the Schnakenberg model is widely studied as a substrate depletion Turing system whilst the Gierer-Meinhardt model is a typical activator-inhibitor Turing system. In both of these cases we consider only one-dimensional domains $\Omega \in (0, 1) $ on which to solve the PDEs and on the boundaries each of the species have no-flux conditions.

\subsubsection{Schnakenberg model}\label{schnakenbergmodel}
The parameterised heterogeneous Schnakenberg model we will be using is as follows.
\begin{align}
    \frac{\partial u}{\partial t} & = \nabla^2u + \gamma\left(- uv^2 + \beta(x)\right) \label{Ch2:Krause:ex1},\\
    \frac{\partial v}{\partial t} & = d \nabla^2v + \gamma\left(uv^2 - v + \eta(x)\right). \label{Ch2:Krause:ex2}
\end{align}
Here $d$ represents the relative diffusion of the activator $v$ compared to that of the substrate $u$ whilst $\beta$ and $\eta$ are spatially dependent production rates. We will focus on a particular form of $\beta$ and $\eta$ in which we parameterise the scale for both the amplitude and frequency of the production heterogeneity 
\begin{align}\label{conditionsSchnaken}
    \beta(x) &= \beta_0 \left(1 + \theta \cos (n \pi x)\right),\\
    \eta(x) &= 1 - \beta(x).
\end{align}
In this way, at each position a combined dimensionless activator/substrate production of 1 is assumed. The parameter $0\leq \beta_0 \leq 1$ describes the average proportion of this production specific to the substrate and the parameter $0\leq \theta \leq 1$ describes the degree of redistribution of the relative production into $n$ periods of peaks and troughs on the domain $\Omega$.
\subsubsection{Gierer-Meinhardt model}\label{GMmodel}
The parameterised heterogeneous Gierer-Meinhardt model is given as follows.
\begin{align}
    \frac{\partial u}{\partial t} & = \nabla^2u + \gamma\left(\frac{u^2}{v}- bu + a(x) \right) \label{Ch2:GM:ex1}\\
    \frac{\partial v}{\partial t} & = d \nabla^2v + \gamma\left(u^2 - v\right), \label{Ch2:GM:ex2}
\end{align}
This model is controlled by the heterogeneous production rate $a(x)$ of the activator $u$. 
We will use a periodic heterogeneity of the form
\begin{align*}\label{conditionsGM}
    a(x) &= a_0 \left(1 +\theta \cos (n\pi x)\right),
\end{align*}
where \(a_0 \in \mathbb{R}\) is the average production rate.

\subsection{Numerical methods}

To generate numerical results we use the numerical continuation method presented by Uecker \cite{Uecker2021} to find solutions of Equations (\ref{MEQ1}) and (\ref{MEQ1BC}) and by starting at $\base{0}$ we find base states for the heterogeneous problem. We begin with the statement that $\Phi(\mathbf{u},\mathbf{x},\mathbf{\bar{x}};\theta) = 0$ ($\mathbf{u}$ must be a solution to Equations (\ref{MEQ1}) and (\ref{MEQ1BC})). By differentiating with respect to $\theta$,
$$ 0 = \frac{\partial \Phi}{\partial \mathbf{u}} \frac{\partial \mathbf{u}}{\partial \theta } + \frac{\partial \Phi}{\partial \theta}.$$ So long as $\frac{\partial \Phi}{\partial \mathbf{u}}$ is nonsingular then $\partial_\theta \mathbf{u}$ can be estimated. As such, finding the base states (and other steady states of the reaction diffusion system) can easily be found using by starting at $\theta = 0$ and incrementing up $\theta$ using a forward Euler approach
\begin{align}
\mathbf{u}_{\theta + \Delta \theta} 
&= \mathbf{u}_\theta + \frac{\partial \mathbf{u}_\theta}{\partial \theta} \Delta \theta, \\ 
&= \mathbf{u}_\theta - \left[ \frac{\partial \Phi_\theta}{\partial \mathbf{u}}\right]^{-1}\frac{\partial \Phi_\theta }{\partial \theta} \Delta \theta \label{updateintheta}
\end{align}
where subscripts indicated the value of $\theta$. The solution generated by Equation (\ref{updateintheta}) is then corrected to reduce error. This is done by setting $\mathbf{u}_{\theta + \Delta \theta} $ as the initial seed of a Newton solver for the problem $\Phi = 0$. We did not find it necessary to use more advanced techniques in increasing $\theta$. 

It is possible to skip the approximate update Equation (\ref{updateintheta}) and simply use a nonlinear solver on $\Phi = 0$ in the vicinity of $\mathbf{u}_\theta$. This is, however, not a good idea since it significantly increases computational time in the nonlinear solver and can sometimes even result in the nonlinear solver finding instead a different steady state solution (of which there may be many). In any case, we make use of the \texttt{pde2path} package which implements this routine. 

Finally, \texttt{pde2path} determines stability by looking at the sign of the largest real component of the eigenvalues of the LHS of the PDE.

In the next section we explore numerical results which give insight into the behaviour of Turing systems with heterogeneous production rates. We first look at the characteristic behaviour of base states (Section \ref{3.1}). Noting that base states often terminate for a sufficiently large value of $\theta$ with a fold bifurcation,  it is clear that for some problems if a heterogeneity is large enough a base state is not defined using our definition. We therefore have a more thorough investigation into what determines if a base state exists or not; what determines how large $\theta$ can be before a fold bifurcation is reached (Section \ref{3.2}). Lastly, how heterogeneous production can affect critical domain lengths required for Turing patterning (Section \ref{3.3}).

\section{Numerical results and discussion}

\subsection{Continuation of steady states}\label{3.1}
The first numerical results illustrate the behaviour of a Schnakenberg Turing system described in Section \ref{schnakenbergmodel} as the heterogeneous production term
is increased in amplitude by tracing the base state and patterned states through numerical continuation of the amplitude parameter $\theta$.
We will first look at some example cases to illustrate the types of branches that can be found.
For all the following results we will use the following parameters; $d=1/40$, $\beta_0 = 0.8$ and $n=1$ unless otherwise stated. Later we will show results for the Gierer-Meinhardt model of Section \ref{GMmodel} where we will use the default parameters $d=20$, $b=1$ and $a_0 = 0.1$ unless otherwise stated. When $\theta =0$ these parameters are known to give a Turing instability in the base state. 
The parameter \(\gamma\) which encodes for the domain length, amongst other things, will be varied between examples to show how the base state behaves as it varies.
In order to visualise the steady state solution branches, we will plot the maximum value on the domain of only the variable \(u\) against the parameter \(\theta\).
This metric has been chosen arbitrarily in order to distinguish between solutions.
It is important to remember when interpreting these bifurcation plots that the branches are only a projection of the infinite dimensional function space onto a single scalar value for plotting purposes.
Importantly, this means that when branches intersect at non-smooth intersections, it is not possible that this is a continuation. Instead, at the point of intersection each branch corresponds to completely unrelated functions (other than the fact that they share a common maximal value of $u$).

In many cases, we observe that there the continuation in $\theta$ can generate base states indefinitely. We can also observe two main bifurcation events on the branch containing the base state. The first of these is a fold at which the base state and the stable patterned state emerge. The second is an example of a fold, terminating the base state, but where the Turing patterned state never bifurcates from base state (they are, instead, perfectly disconnected). By saying `patterned state' we are implying that there is a branch corresponding to a non-homogeneous but also stable steady state (indicated in blue in each figure). Finally, we demonstrate some exotic behaviour of the steady states under some conditions.

\subsubsection*{Base state with no limitation}
\begin{figure}
    \centering
    \includegraphics[width=\textwidth]{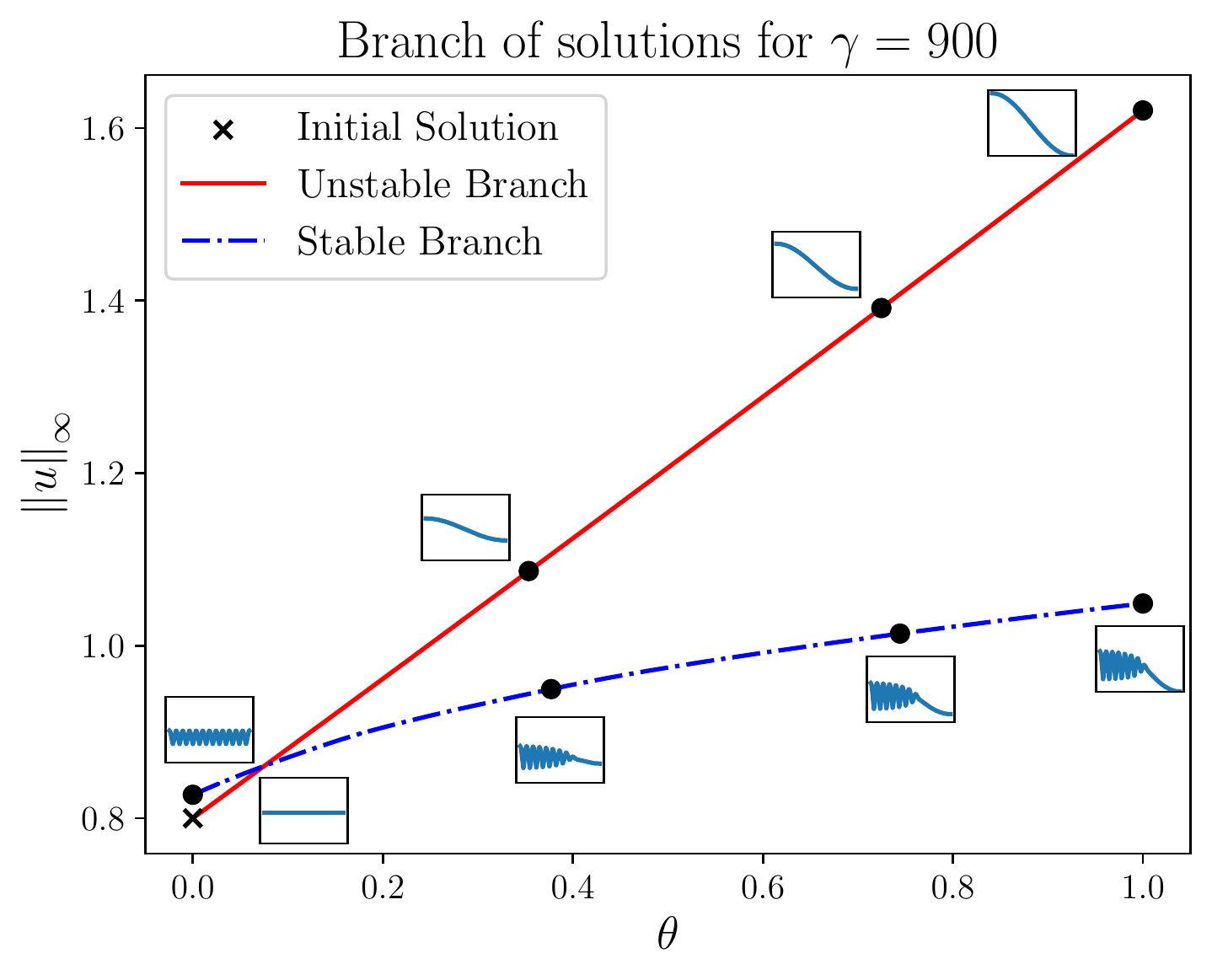}
    \caption{Schnakenberg system bifurcation diagram for growing heterogeneity $\theta\in[0,1]$. Parameters used are characteristic of \textit{large} domains relative to Turing pattern wavelength ($\gamma = 900$) with also $\beta_0 = 0.8$ and  $d=1/40$. When $\theta = 0$, the system solves a classical Turing system where the base state is homogeneous and indicated with an $\times$. As the heterogeneity $\theta$ grows, so does the base state. A number of examples of the spatial distribution of $u$ along the (red) unstable base state $\base{\theta}$ is displayed. In this case, the base state is allowed to grow continuously without a fold. On the other hand, a (blue) stable Turing `patterned' state branch is also shown with some displayed distributions of $u$. This is found by solving the full reaction-diffusion equation at $\theta=0$ and applying the numerical continuation.
    }
    \label{fig:fold4}
\end{figure}

In the most simple case, starting with $\base{0}$ and growing the heterogeneous term by increasing $\theta$ in Section \ref{schnakenbergmodel}, no folds were found in increasing $\theta$ from 0 to 1. It is important to note that this does not mean that the base states will extend for an arbitrarily large $\theta$. For the Schnakenberg system in Section \ref{schnakenbergmodel}, we find that this often occurs for large $\gamma$ and in Fig. \ref{fig:fold4} use the value of $\gamma = 900$. This corresponds with a very large domain in relation to the expected wavelength of any Turing patterns.
Our value of $\gamma$ corresponds to a value of $\epsilon \approx 1.1 \times 10^{-3}$ in the paper by Krause \textit{et al.} \cite{Krause2020}. We find that in this case the base state exists by numerical continuation and furthermore that it is approximately equal to the steady state where diffusion is ignored as small which is trivial because it is clear from Equations (\ref{MEQ1}) and (\ref{MEQ1BC}) that unless $\theta$ large on the order of $\gamma$, for large $\gamma$, we simply have to leading order that $\base{\theta}$ solves $F(\mathbf{u})+\theta G(\mathbf{u},x) = 0$.

\subsubsection*{Base state fold connected to a patterned state}
\begin{figure}
    \centering
    \includegraphics[width=\textwidth]{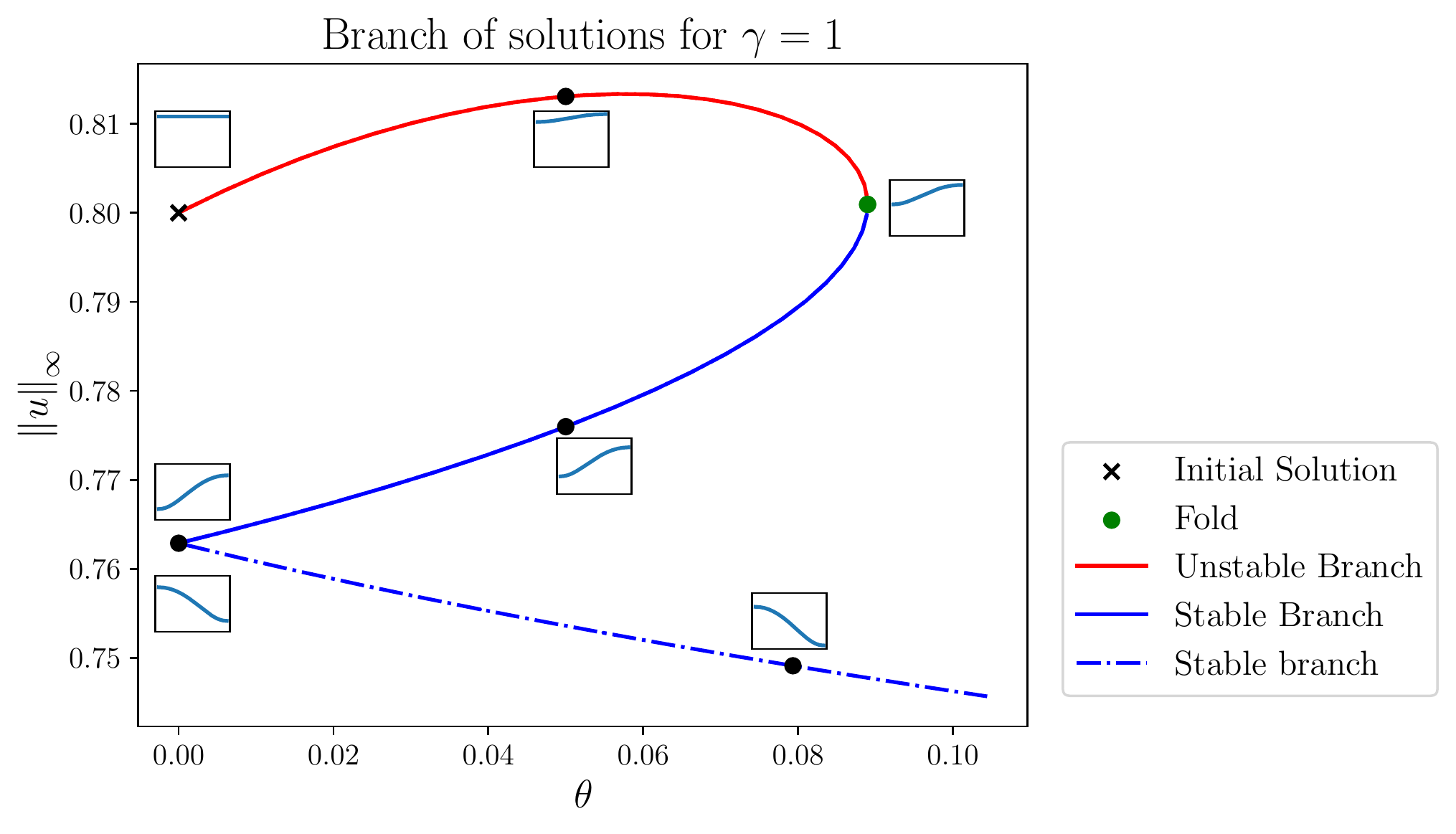}
    \caption{Schnakenberg system bifurcation diagram for growing heterogeneity $\theta\in[0,1]$. Parameters used are characteristic of \textit{small} domains relative to Turing pattern wavelength ($\gamma = 1$) with also $\beta_0 = 0.8$ and  $d=1/40$. When $\theta = 0$, the system solves a classical Turing system where the base state is homogeneous and indicated with an $\times$. As the heterogeneity $\theta$ grows, so does the base state. A number of examples of the spatial distribution of $u$ along the (red) unstable base state $\base{\theta}$ is displayed. In this case, the base state merges with the stable patterned state at around $\theta = 0.09$. The blue branches are stable patterned states but only the solid branch can be obtained by continuing through the fold. The dot-dash branch can be found through continuation of a fold in the base state if \textit{decreasing} $\theta$ from the $\theta = 0$ base state.
   }
    \label{fig:fold1}
\end{figure}

We observe different behaviour in the base state for non-large $\gamma$. If $\gamma$ is small, but not too small as to not observe Turing patterns in the homogeneous Schnakenberg system (due to the domain size being less than the necessary critical domain length), then we observe a critical fold in the base state solution. 
In Fig. \ref{fig:fold1}, we use the value of $\gamma = 1$. When $\theta = 0$, this corresponds to the case where there is just one unstable wavenumber corrsponding to a Turing pattern with just a half period on the full domain.  
In this case, the branch for a patterned state merges with the branch of the base state, undergoing a fold bifurcation as seen in Fig. \ref{fig:fold1}.
This means that the base state becomes closer and closer to a patterned state until both states are indistinguishable from each other at the fold bifurcation. For heterogeneities with an amplitude $\theta$ beyond this fold (shown with a green dot in Fig. \ref{fig:fold1}), we are unable to objectively define a suitable base state and therefore it becomes ambiguous as to whether or not a `Turing' pattern is observed in the solution of the reaction-diffusion problem. Indeed, whilst a steady state solution to the reaction-diffusion equation is expected beyond the fold, we do not know where this solution is by numerical continuation from $\theta=0$ without significant work. That is, there are other missing branches here and it remains unclear if any of these are reasonable candidates to be defined as a `base state' at this stage and further work here is needed. In Fig. \ref{fig:fold1}, you can see the stable patterned state but also an unstable patterned state. For $\theta=0$ there are at least two patterned states. You can see these states in the bifurcation diagram as mirrored functions. Interestingly, if the heterogeneity is inverted in sign ($\theta\in[-1,0]$), continuation shows a mirror image of the bifurcation diagram in Fig. \ref{fig:fold1}.

\subsubsection*{Base state fold not connected to a patterned state}
\begin{figure}
    \centering
    \includegraphics[width=\textwidth]{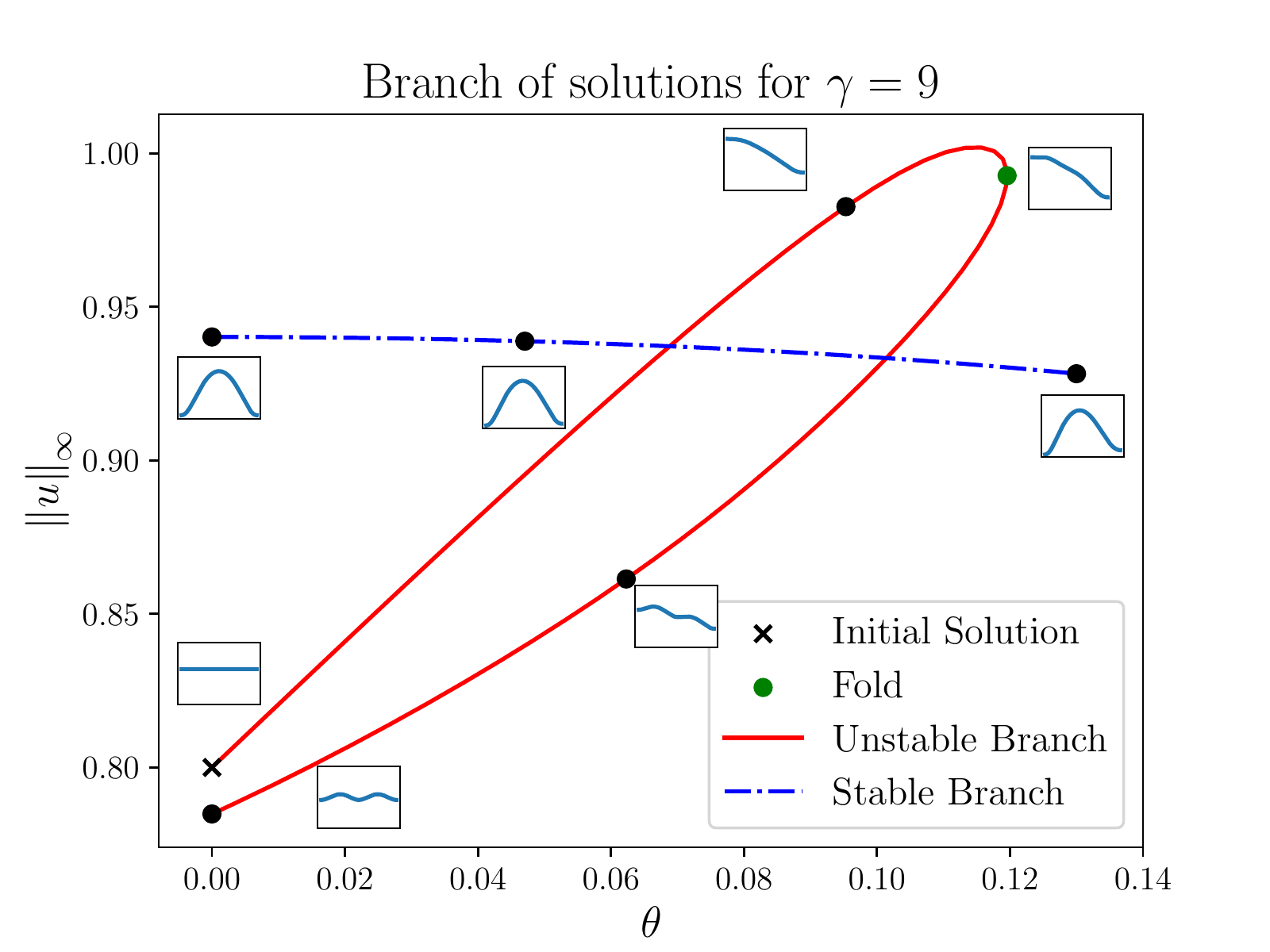}
    \caption{Schnakenberg system bifurcation diagram for growing heterogeneity $\theta\in[0,1]$. Parameters used are characteristic of \textit{intermediate} domains relative to Turing pattern wavelength ($\gamma = 9$) with also $\beta_0 = 0.8$ and  $d=1/40$. When $\theta = 0$, the system solves a classical Turing system where the base state is homogeneous and indicated with an $\times$. As the heterogeneity $\theta$ grows, so does the base state. A number of examples of the spatial distribution of $u$ along the (red) unstable base state $\base{\theta}$ is displayed. In this case, the base state merges with an unstable heterogeneous steady state at around $\theta = 0.12$. The blue branch is a stable patterned state but the dot-dash nature of this branch indicates that it is not obtained by continuation past a fold from the steady state but instead by solving the reaction-diffusion equation with $\theta = 0$ until steady state and using continuation from there.}
    \label{fig:foldUnstable}
\end{figure}

In intermediate values of $\gamma$, more curious behaviour is possible. This is in part because these values permit multi-wavelength heterogeneous steady states. In Fig. \ref{fig:foldUnstable} we now display the bifurcation diagram for \(\gamma = 9\) (analogous to a domain length increase of three-fold on the example in Fig. \ref{fig:fold1}).
The key observation in Fig. \ref{fig:foldUnstable} is that whilst the base state branch also undergoes a fold bifurcation,
 the solution branch with which it merges is an unstable heterogeneous steady state (not a stable pattern).
This illustrates that the base state branch can merge with another branch which is not a branch of patterned states. In considering Fig. \ref{fig:fold4} where the base state seemingly continues indefinitely without folds, it is possible that a fold is present in a similar way to how it appears in Fig. \ref{fig:foldUnstable} but at sufficiently large values of $\theta$. If this is the case, our observations might suggest that as $\gamma$ gets very large, so to does the values of $\theta$ where base state folding first occurs. 

\subsubsection*{Exotic behavior}
\begin{figure}
    \centering
    \includegraphics[width=\textwidth]{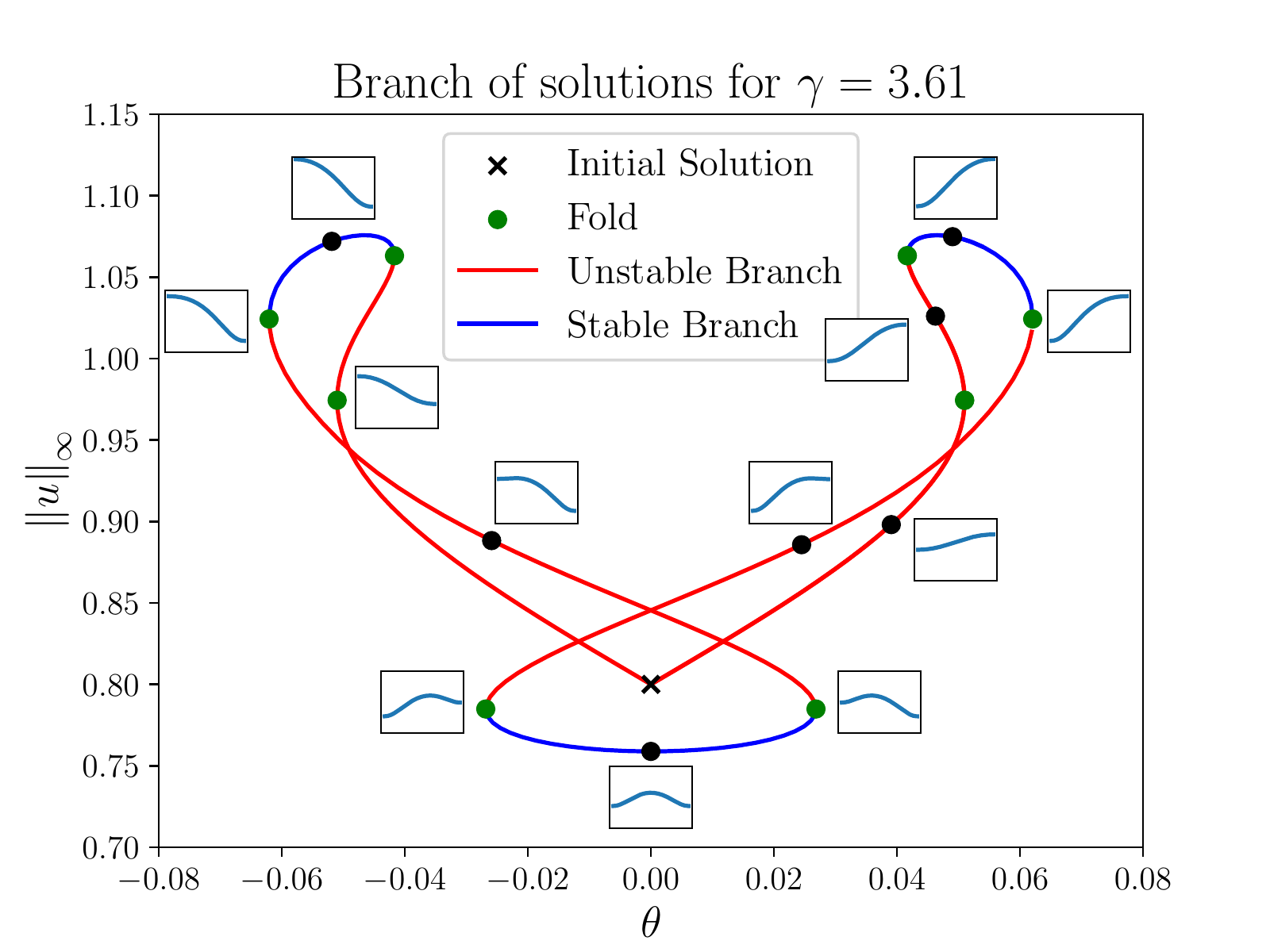}
    \caption{Schnakenberg system bifurcation diagram for growing heterogeneity $\theta\in[-1,1]$. Parameters used are characteristic of \textit{narrowly defined} domains relative to Turing pattern wavelength ($\gamma = 3.61$) with also $\beta_0 = 0.8$ and  $d=1/40$. When $\theta = 0$, the system solves a classical Turing system where the base state is homogeneous and indicated with an $\times$. As the heterogeneity $\theta$ grows, so does the base state. A number of examples of the spatial distribution of $u$ along the (red) unstable base state $\base{\theta}$ is displayed. Note that here the base state would only be defined between approximately -0.05 and 0.05. By continuing through each fold, we end up back at $\base{0}$. Interestingly, this closed loop contains three different patterned branches (blue) but not a patterned branch on approximately $\pm(0.03,0.04)$. It is expected that the patterned state obtained by solving the reaction-diffusion equation in this regime is not connected here.}
    \label{fig:fold2}
\end{figure}

While the previous examples show two branches originating at \(\theta = 0\) converging, this does not capture all possibilities.
In a more bizarre scenario, we can consider the case where \(\gamma = 3.61\).
As shown in Fig. \ref{fig:fold2}, the system undergoes many folds before merging with another solution branch which contains \(\theta = 0\). Furthermore, there are stable steady states which are only present for a discrete range of \(\theta\) values. To demonstrate the behaviour and the way it closes itself, it was necessary to continue in both the positive and negative $\theta$ direction from $\base{0}$.

\subsection{Base state existence} \label{3.2}

\begin{figure}

    
    \centering
    \includegraphics[width=\textwidth]{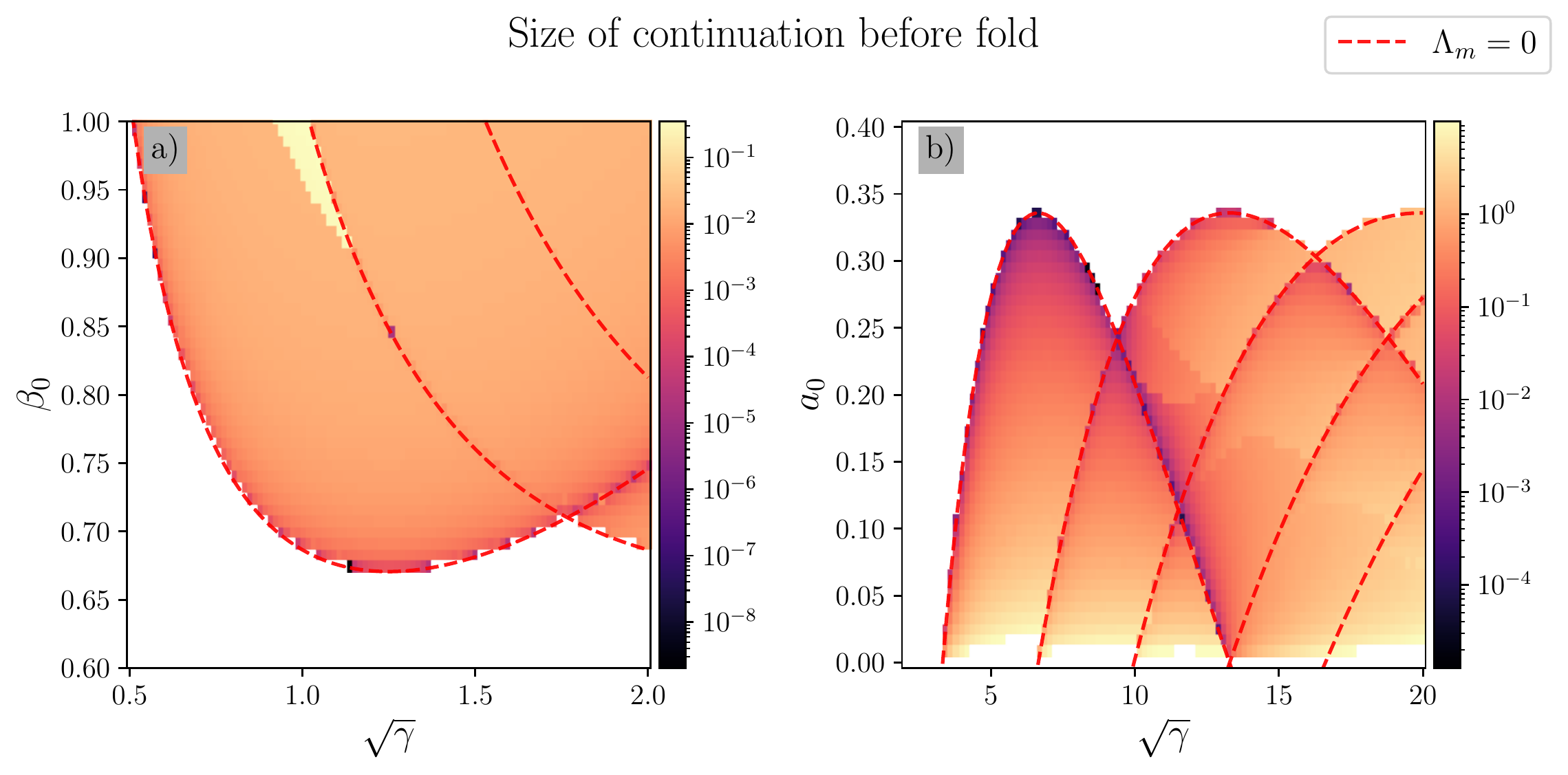}
    \caption{Size of continuation before a fold $\theta^+$ for (a) the Schnakenberg model and (b) the Gierer-Meinhardt model as $\gamma$ is varied along with (a) $\beta_0$ and (b) $a_0$, respectively. The size of the continuation is presented in color on the log scale. All of these results are given for $n=1$ in the heterogeneous term in the respective models. Red curves are drawn on the figures to correspond with \(\Lambda_m = \mathrm{max}_j \Re\left(\lambda_{j}(\mathsf{A}_{m}) \right) = 0\) for $m=1,2,3$ (for curves left to right on both subfigures) where $\lambda_{j}(\mathsf{A}_{m})$ are eigenvalues defined in Ssection \ref{BaseStates}. The background color of white indicates that no fold was found for these parameter sets and $\theta$ was allowed to grow to 1.
        }
    \label{fig:foldScan1}
\end{figure}

\begin{figure}
    %

    
    \centering
    \includegraphics[width=\textwidth]{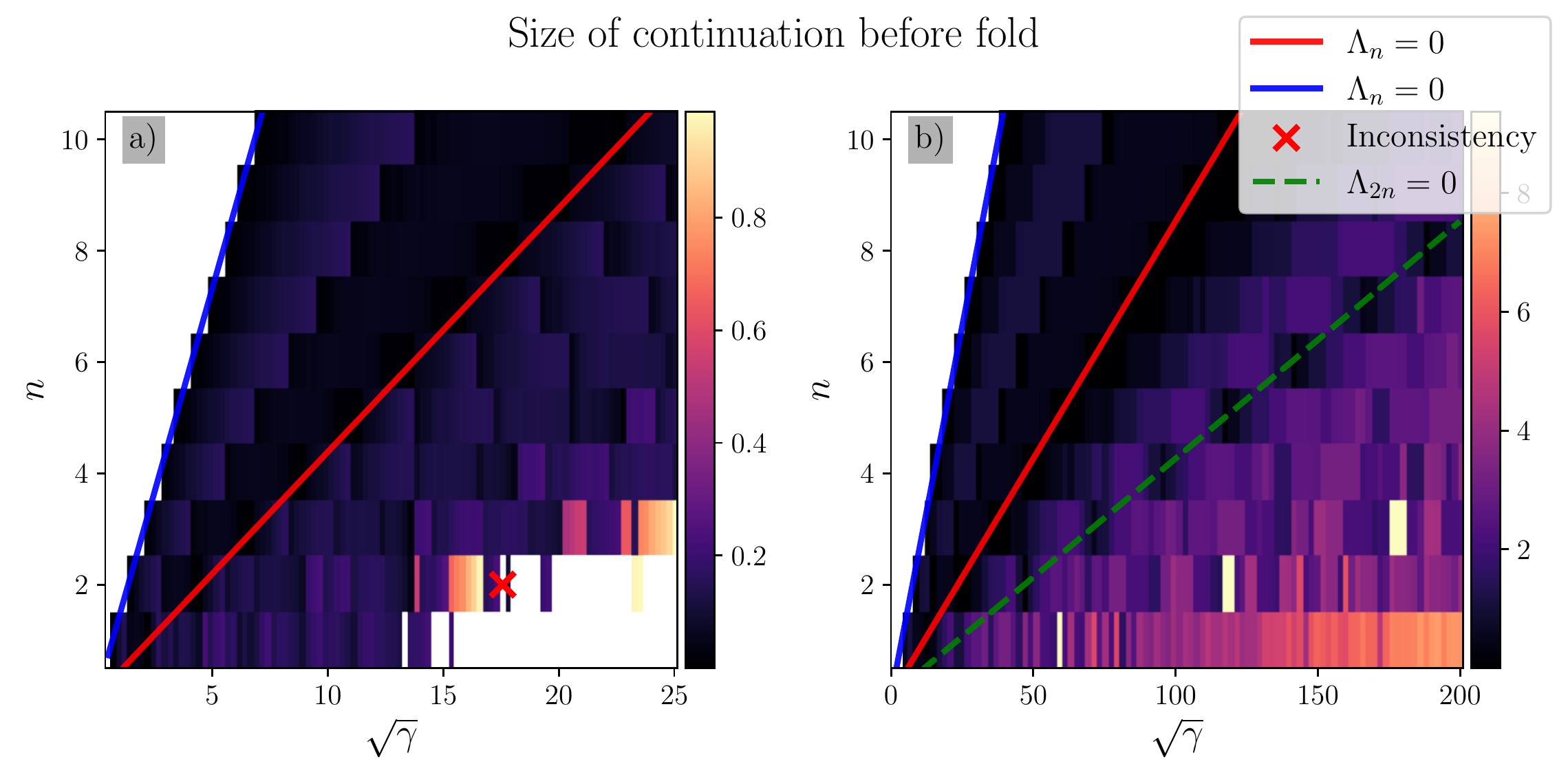}
    \caption{Size of continuation before a fold $\theta^+$ for (a) the Schnakenberg model and (b) the Gierer-Meinhardt model as $\gamma$ and $n$ is varied for each model. The size of the continuation is presented in color. Setting (a) $\beta_0 = 0.8$ and (b) $a_0 = 0.1$ in each model respectively, \(\Lambda_n = \mathrm{max}_j \Re\left(\lambda_{j}(\mathsf{A}_{n}) \right) = 0\) where $\lambda_{j}(\mathsf{A}_{m})$ are eigenvalues defined in Section \ref{BaseStates} has two solutions. The solution with smallest $\gamma$ is shown on the blue line and the other is shown on the red line. The background color of white indicates that no fold was found for these parameter sets and $\theta$ was allowed to grow to 1. In (b) the green dashed line is an overlay of the red line with half of the value of $n$ for each $\gamma$. This curve surprisingly traces a pattern of small $\theta^+$. In (a) a red $\times$ indicates a continuation that runs into numerical difficulties.}
    \label{fig:foldScan2}
\end{figure}

In order to have a discussion about Turing patterns, it is important for a base state to exist. It is therefore critical to explore what determines $\theta^+$, the maximum size that $\theta$ can take before a critical point such as a fold is encountered. To accomplish this we performed parameter scans on both the Schnakenberg and Gierer-Meinhardt model from Sections \ref{schnakenbergmodel} and \ref{GMmodel}. Our immediate observation from doing these scans is that fold bifurcations are very common. In particular, we observed more folds when the spatially-dependent source term $\mathbf{G}(\mathbf{u},x)$ varies explicitly in space with frequencies similar to that of unstable eigenvectors in the dispersion relation.

In Fig. \ref{fig:foldScan1} we look at $\theta^+$ for the Schnakenberg model (a) and the Gierer-Meinhardt model (b). In Fig. \ref{fig:foldScan1} (a) we plot $\theta^+$ as the scale parameter $\gamma$ and the parameter $\beta_0$ in the Schnakenberg model are varied, whilst in (b) we instead vary the parameter $\alpha_0$ in the Gierer-Meinhardt model. In both cases, we have plotted, in red, the curves that relate to eigenvalues \(\Lambda_m = \mathrm{max}_j \Re\left(\lambda_{j}(\mathsf{A}_{m}) \right) = 0\) for $m=1,2,3$ (for curves left to right). We note that in our test problems we do not have strictly imaginary eigenvalues so along these curves $\mathbf{\bar{J}}_0$ is singular and we expect that $\theta^+$ is not finite.  For each constant $\beta_0$ (or $\alpha_0$) we see that $\Lambda_m = 0$ at most twice because solving $\Lambda_m = 0$ requires solving a quadratic. Between the two values, we find that $\Lambda_m > 0$ and thus the $m$th mode of the homogeneous problem is unstable. On these curves, $\mathbf{\bar{J}}_0$ is singular. As previously established, we expect on these curves that continuation is not possible. In the region shown in white, we found no upper bound in $\theta^+$. This region also corresponds to the subset of the parameter space where the associated homogeneous system is devoid of Turing patterning. The red curves furthest to the left correspond to $m=1$ (corresponding to the onset of Turing instability in the eigenfunction $\cos(\pi x)$ at $\theta =0$). Note that our growing heterogeneity is also of this form ($n=1$) $\cos(n \pi x )$ (see Equations (\ref{conditionsSchnaken}) and (\ref{conditionsGM})). We find that because of this a fold is very quick to form in the numerical continuation near the red curve corresponding to $m=1$ but not near the onset of instability for the higher modes. Small $\theta^+$ is shown by darker colors in the plot.

To investigate specifically if small $\theta^+$ is associated with $m=1$ because $n=1$ we varied $n$ in the Schnakenberg model from 1 to 10. In Fig. \ref{fig:foldScan2}, for each $n$, holding $\beta_0 = 0.8$ (a) and $\alpha_0 = 0.1$ (b) we plot the size of the continuation $\theta^+$ as $\gamma$ is increased. We indicate the minimum value of $\gamma$ (blue line) and the maximum value of $\gamma$ (red line) for which $\Lambda_n = 0$. That is, for $n=1$ the blue and red curves correspond to the first and second intersection of $\beta_0 = 0.8$ (a) and $\alpha_0 = 0.1$ (b) with the respective red curves in Fig. \ref{fig:foldScan1}. We see for each $n$, the size of $\theta^+$ is very small at both zeros of $\Lambda_n$. What is also surprising, if $n$ is larger than 1 if $\gamma$ is smaller than that required to make the $n$th mode unstable in the homogeneous problem, the continuation did not fold. That is, we may have a Turing instability in the homogeneous problem because of an instability in the $m=1$ mode but if the heterogeneity has a higher spatial frequency, say $n=2$, the base state may not encounter a fold readily. As the scale parameter $\gamma$ is increased beyond the the red line, we find what appears to be noise in $\theta^+$ but within this noise appears to be patterns. Looking specifically at the Gierer-Meinhardt model in Fig. \ref{fig:foldScan2} (b) we see small $\theta^+$ near the value of the maximum $\gamma$ for which $\Lambda_{2n} = 0$. We have indicated that this is case by tracing the green dashed line over the expanse of small $\theta^+$. You can also see this effect in Fig. \ref{fig:foldScan1} (a) for $n=1$ by looking at the left branch of the $m=2$ red curve and seeing a noticeable dark shade. As \({\gamma}\) increases, the magnitude that $\theta$  can be continued before reaching a fold tends to increase, before not reaching a fold at all.
However, numerical instabilities are prevalent in this region, as shown specifically by the red $\times$ in Fig. \ref{fig:foldScan2} (a), so the accuracy of these results remains questionable. We shall look specifically at the continuation described by this red $\times$ in the next section.
The numerical results seem to become more accurate as the spatial grid becomes finer, and the maximum step size in $\theta$ becomes smaller.
Due to the computational cost of producing parameter scan results, the accuracy of the results is here limited.

\subsubsection*{Numerical Issues}
\begin{figure}
    \centering
    \includegraphics[width=0.95\textwidth]{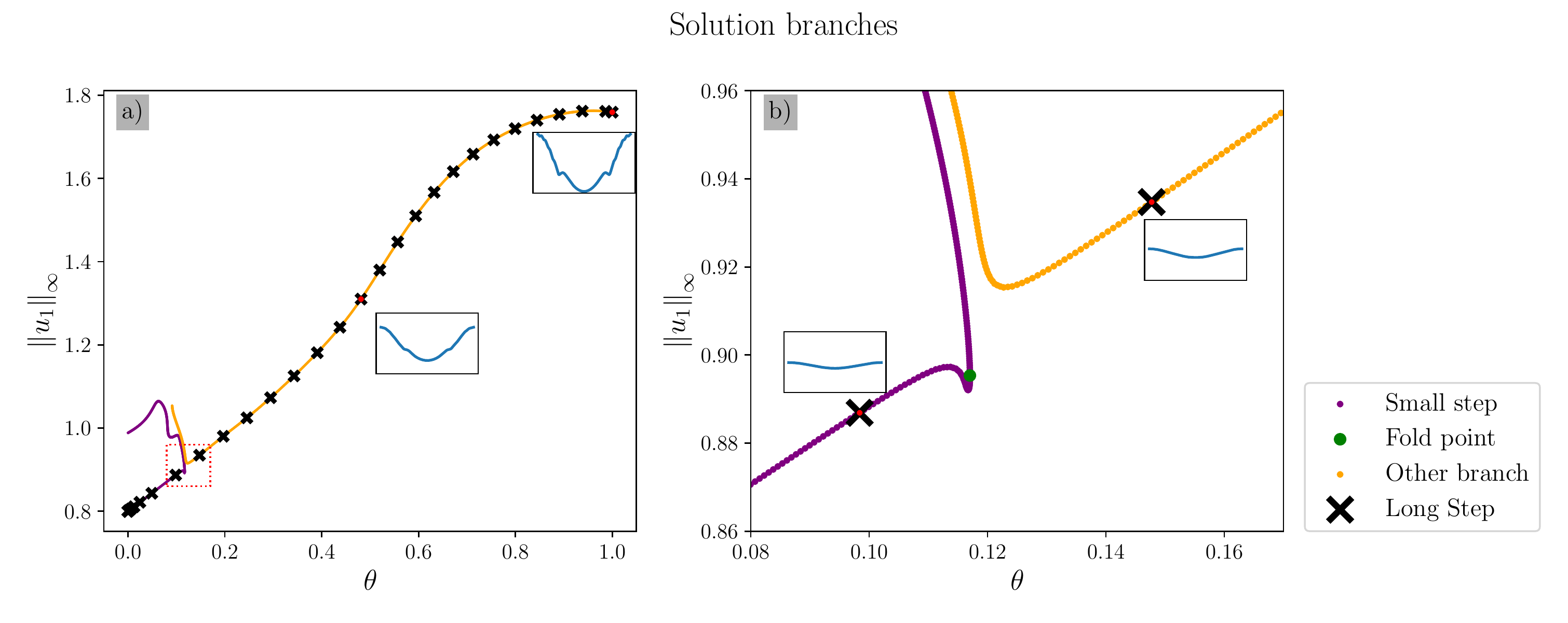}
    \caption{Plot of branches for the numerically inconsistent case highlighted in Fig. \ref{fig:foldScan2} (a) with varying maximum step size. In purple, the base state branch and continuation through the fold point (green dot) with very small step sizes is shown. In yellow, a different branch is shown and the $\times$ symbols show the updates in the continuation algorithm if the step size is too coarse. Plot (a) shows the full bifurcation diagram whilst plot (b) displays a zoomed version of the region enclosed in the red box to show detail near the fold point.}
    \label{schScan2Inconsistency}
\end{figure}

The inconsistent numerical issue that occurs occasionally in our parameter sweeping experiments in the previous section are investigated here. In particular, we investigate the red $\times$ continuation in Fig. \ref{fig:foldScan2} (a).
In this continuation a maximum step size of \(10^{-1}\) was used. This is a relatively large step size, but since the~\texttt{pde2path} package adaptively adjusts the step size as needed, it can usually make out the finer details without much increase in computational cost.
However, in this case, the larger step size causes the solution to jump from one branch to another. This can be seen in bifurcation Fig. \ref{schScan2Inconsistency}, where for a small step size, a fold is encountered early in the continuation, but for a large step size, the continuation jumps to a different branch.
Clearly the results in this region are unreliable.
It is not clear how small the step size must be made in order to avoid this occurring. It does raise an interesting question though. In this example, it is pretty clear that the (yellow) branch that the coarse numerical algorithm found does not technically satisfy the numerical continuation criteria for a base state. That being said, looking at the distributions on either side of the singularity, it is possible that the yellow branch perhaps \textit{should} be considered a base state. It remains unclear if such a suitable branch can be found in for other cases. However, this case hints at the possibility that there may be a better definition for a base state  than the one presented in this manuscript (one which can potentially always describe a unique state for all problems).

\subsection{Critical domain length} \label{3.3}
The extension of the Turing instability to spatially-dependent RD systems allows us to distinguish between patterned states and the base states.
Previously these solution states were often indistinguishable.
This meant that analysing certain phenomena, such as the critical domain length, was very challenging or impossible.
Now that the Turing instability has a spatially-dependent analogue, we can study such phenomena.
As a proof of concept, we will study how the critical domain length changes as the size of the heterogeneity in a spatially-dependent RD system increases.
The critical domain length has important physical implications, especially in developmental scenarios.
In a scenario where the domain is slowly growing, Turing patterns will arise only if the size of the domain is above the critical domain length.
Therefore, assessing the impact of a spatially-dependent term on the critical domain length could have key implications for these developmental scenarios. 
We will attempt to investigate the change in the critical domain length with respect to the size of the heterogeneity for two different reaction terms.

The critical domain length is encoded in a critical \(\gamma\) value which we will call \(\gamma_c\).
Denote \(\gamma_{c,0} \in \mathbb{R}^{+}\) as the critical \(\gamma\) value for the classical RD system, and \(\gamma_{c,\theta} \in \mathbb{R}^{+}\) as the critical \(\gamma\) value for the heterogeneous RD system with parameter $\theta$.
Further, define \(L_{c,0} \coloneqq \sqrt{\gamma_{c,0}}\), \(L_{c,\theta} \coloneqq \sqrt{\gamma_{c,\theta}}\) as the respective critical domain lengths.
Here we are accepting \(L_c = \sqrt{\gamma_c}\) to be a non-dimensional equivalent of the critical domain length. 

The value of \(\gamma_{c,\theta}\) is defined by largest \(\gamma\) such that the base state of~Equations (\ref{Ch1:parameterised}) and (\ref{Ch1:parameterisedBC}) is stable for all \(\gamma < \gamma_{c,\theta}\), but exhibits Turing instabilities for some \(\gamma > \gamma_{c,\theta}\). 
It is infeasible to check all \(\gamma\) values less than some candidate value for \(\gamma_{c,\theta}\).
Instead, we can rely on the fact that when \(\gamma = \gamma_{c,0}\), $\Lambda_m = 0$ which can be calculated exactly for both the Schnakenberg model and Gierer-Meinhardt model.

Instead of parameterising the base state branch with the size of the heterogeneity \(\theta\) only, we will also parameterise with respect to \(\gamma\).
In doing so, we are assuming that a path independence result holds.
That is, the base state solution for some \(\gamma_0 > 0\) can be found by first finding the base state solution for another \(\gamma_1 > 0\), and then continuing from that base state solution with respect to \(\gamma\) to find the solution at \(\gamma_1\). 
Initially we will use \(\gamma = \gamma_{c,0}\) to perform the continuation, as this is known exactly and we will assume that this is close to \(\gamma_{c,\theta}\).
After finding a base state solution with the initial \(\gamma\) value, we perform numerical continuation with respect to \(\gamma\), and continue to increasing or decreasing $\gamma$ until finding \(\gamma_{c,\theta}\) for a given $\theta$.
We reach the critical value $\gamma_{c,\theta}$ when the base state (with respect to $\gamma$ but constant $\theta$) undergoes a change of stability.
If the base state found for \(\gamma = \gamma_{c,0}\) is stable, then we will increase \(\gamma\) in the second stage continuation.
Likewise, we will decrease \(\gamma\) if the base state is unstable. Determining whether a steady state solution is stable can be done using inbuilt methods in~\texttt{pde2path}~\cite{Uecker2021}.

We are relying on using \(\gamma = \gamma_{c,0}\) as an initial condition for the continuation.
However, based on recent analysis on heterogeneous RD systems, there are points where the system with \(\theta = 0\) is outside of the Turing region, so we still expect to see Turing instabilities for a sufficiently large \(\gamma\)~\cite{Krause2020}.
If the homogeneous system defined by \(\theta=0\) is outside of the Turing region, it is unclear what the initial \(\gamma\) value should be.
A further investigation into resolving a method for finding the critical domain length in this case should be considered.


\begin{figure}
    \centering
    \includegraphics[width=\textwidth]{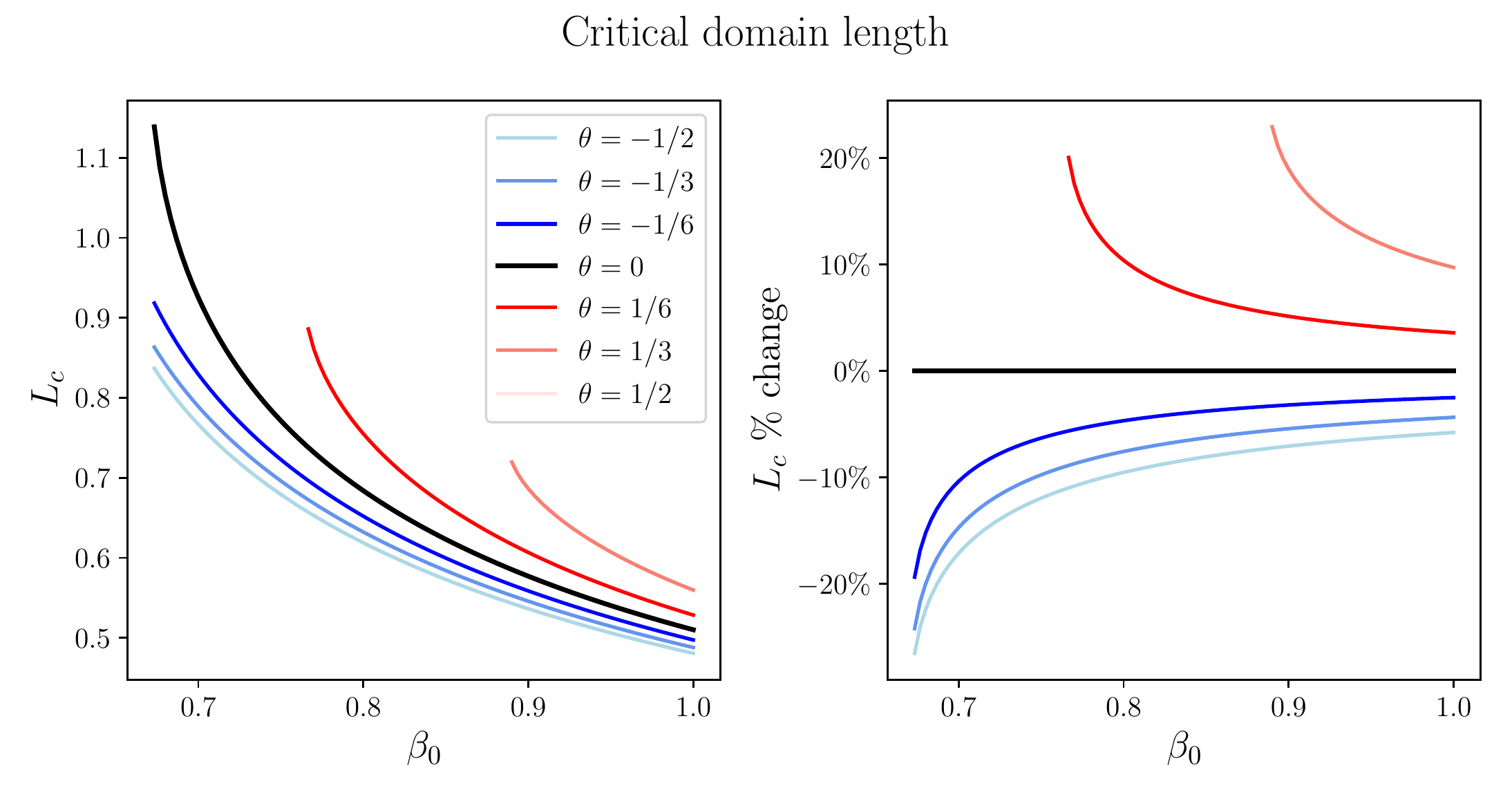}
    \caption{Critical domain lengths \(L_{c,\theta}\) of the Schnakenberg system described in Section \ref{schnakenbergmodel}. The critical domain length is plotted for a range of heterogeneity sizes \(\theta\) as a function of the parameter \(\beta_0\).}
    \label{fig:derivSch}
\end{figure}

\begin{figure}
    \centering
    \includegraphics[width=\textwidth]{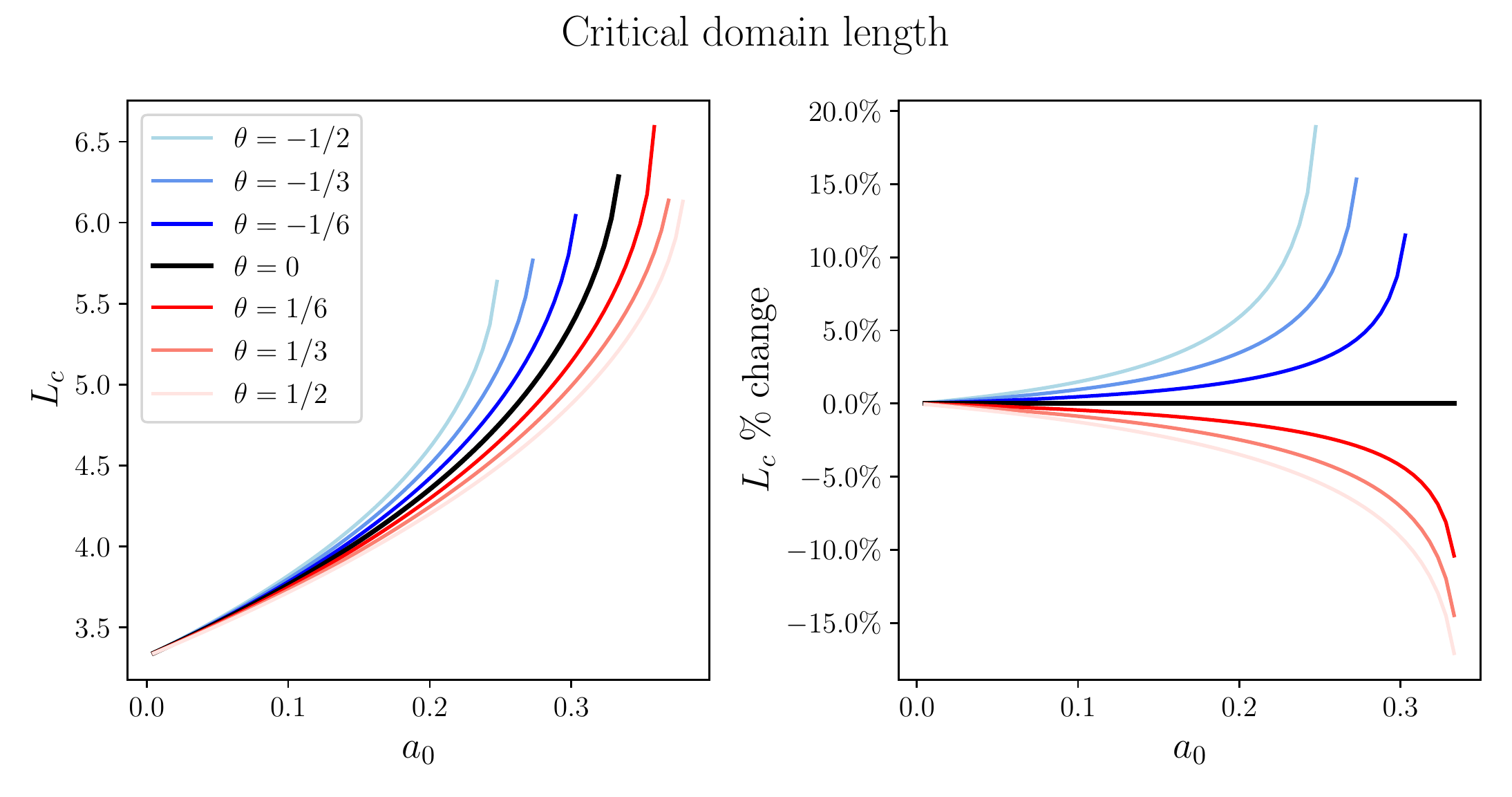}
    \caption{Critical domain lengths \(L_{c,\theta}\) of the Gierer-Meinhardt system described in Section \ref{GMmodel}. The critical domain length is plotted for a range of heterogeneity sizes \(\theta\) as a function of the parameter \(a_0\).}
    \label{fig:derivGM}
\end{figure}

Fig. \ref{fig:derivSch} shows the critical domain length \(L_c\) for the Schnakenberg system for a range of \(\theta\) and \(\beta_0\) values.
The length \(L_c\) appears to be decreasing with respect to \(\beta_0\) and increasing with respect to \(\theta\).
On the other hand, Fig. \ref{fig:derivGM} shows that the critical domain lengths for the Gierer-Meinhardt system appears to have the reverse dependence on the parameter $a_0$.

For a given production rate, if the \(\theta = 0\) is within the Turing region, then we expect to have a critical domain length for every other \(\theta\) value.
This is because the cosine heterogeneity will cause at least one interval of the domain to be within the Turing region locally.
Thus, for sufficiently large \(\gamma\), we expect to see Turing patterns~\cite{Krause2020}.
However, our method for finding the critical domain length in many of these cases fails.
Most notably, the critical domain length could not be found for any \(\beta_0\) value when \(\theta = 1/2\), as seen in~Fig. \ref{fig:derivSch}.
This is potentially because there is a decoupling effect between two intervals which are locally within the Turing region.
\begin{figure}
    \centering
    \includegraphics[width=\textwidth]{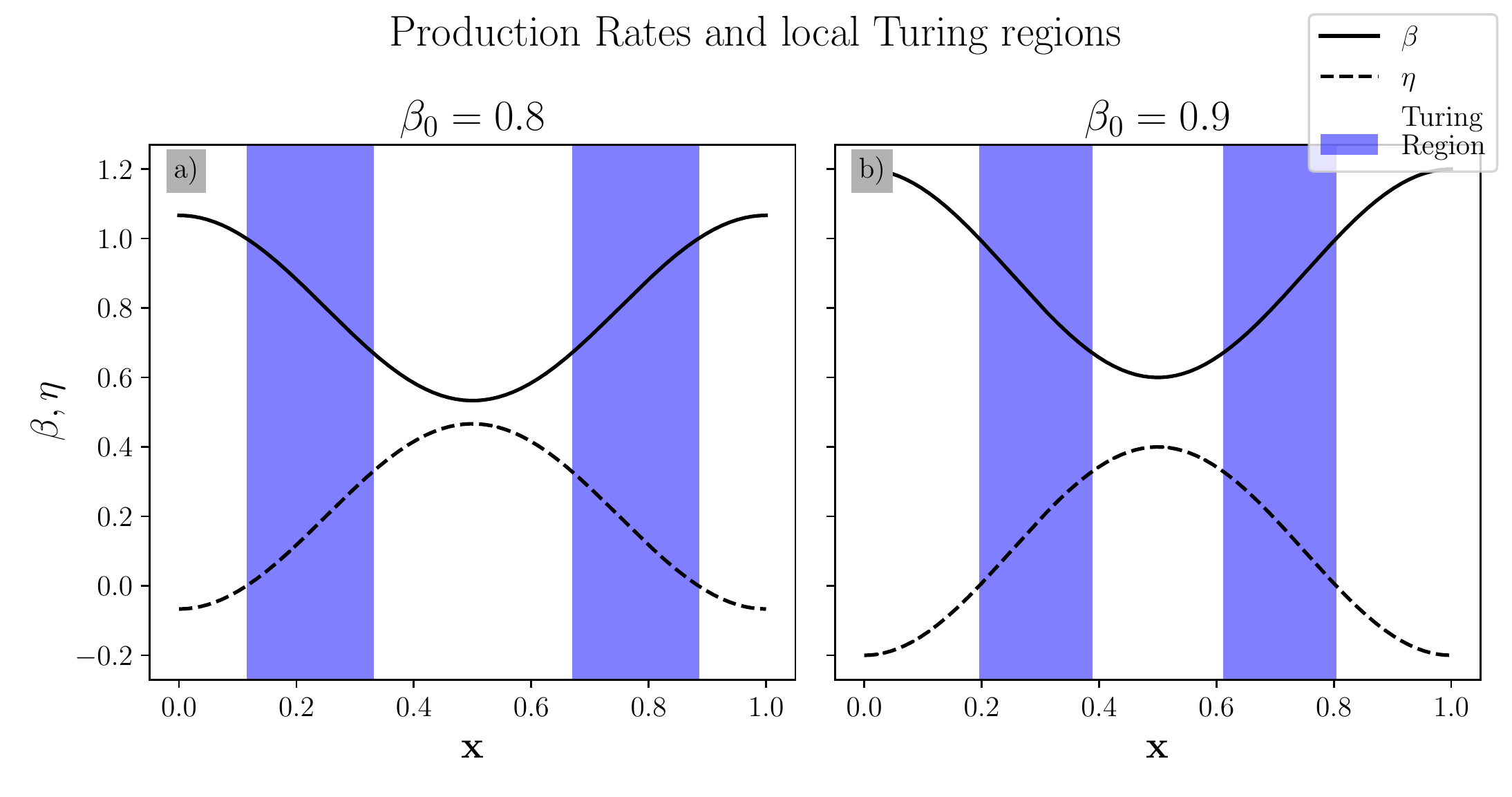}
    \caption{Production rates for the first chemical, \(u \), and the second chemical \(v\), for the Schnakenberg model of Section \ref{schnakenbergmodel}. Plots (a) and (b) describe the model with \(\beta_0 = 0.8\) and \(\beta_0 = 0.9\) respectively. Each figure also shows the regions where the system is \textit{locally} within the classical Turing pattern-generating parameter space. These plots are made for \(\theta = 1/3\), meaning that we found a critical domain length for the system shown in (b), but not in (a). In (a), gap between the regions that are driving the Turing instability in the whole domain are further apart and it is possible that these are effectively decoupled. In this case, we would expect to find a critical domain length but significantly larger (where Turing patterns can be associated with the sub-domains which locally drive Turing patterns).
    }
    \label{fig:critDomainLengthDiscussion}
\end{figure}
Fig. \ref{fig:critDomainLengthDiscussion} shows the regions where the systems with \(\theta = 1/3\) and \(\beta_0 = 0.8, 0.9\) are locally within the Turing region.
As seen in Fig. \ref{fig:derivSch}, a critical domain length could be found for \( \beta_0 = 0.9 \), but not for \(\beta_0 = 0.8\).
Although the Turing regions are larger in the case where \(\beta_0 = 0.8\), the region between the two Turing regions is also larger.
This gap between the Turing regions could have a decoupling effect where, if the two regions are close enough together, they can act as one region for the purposes of forming a Turing instability.
That is, there is enough bleed through from one region to the other to support a Turing pattern, despite having a region where no Turing pattern can be supported in between.
So in this case there would be a critical \(\theta\) value after which \(\gamma\) must be significantly larger before observing Turing instabilities which are local to the respective Turing regions.

\section{Conclusions}
\label{sec:conclusions}
Despite being widely applicable to various problems in science, Turing instabilities in spatially-dependent reaction-diffusion systems have yielded very little attention in the literature.
One of the roadblocks to understanding the behaviour of these systems is the lack of definition for Turing instabilities when the problem depends on the spatial coordinate.
The classical definition relies on the existence of a uniform steady state solution, however no such steady state exists for spatially-dependent problems in general.
In reformulating the definition, the problem arises of distinguishing between patterned states and the base state.
The base state in the classical case is the uniform steady state.
Since the steady state solutions of most spatially-dependent reaction-diffusion system are non-uniform, it is unclear which states we should label as `patterned', and which are labelled as a `base state'.
In order to link the spatially-dependent case with the classical case, we utilise tools from continuation to gradually increase the size of the heterogeneity.
That is, the spatially-dependent term (or heterogeneity), is parameterised such that the heterogeneity vanishes initially, and grows to full amplitude as the introduced parameter increases.
Once at full amplitude, the base case solution to the reaction-diffusion equation is the solution found through continuation, with a full amplitude heterogeneity.
This grounds the spatially-dependent base case to the classical base case, and allows us to distinguish between patterned and non-patterned states.
By defining the base case solution through continuation, this also provides a method for finding the base solution using numerical continuation.

While we have extended the definition of the Turing base state, this does not directly extend the definition of the Turing instability.
Traditionally, a Turing instability requires the base state to be stable to constant perturbations, and unstable overall.
The stability to constant perturbations condition is not relevant with a spatially-dependent base state.
As such, the extension of the first Turing condition is not trivial even after defining the base state.
So we discussed a few possibilities about how this condition could be extended, and the benefits of each possibility. Much more research can be done to analyse the properties of each of these definitions.

After defining the base state for heterogeneous Turing systems, it remains to know whether such base states exist.
We provided a variety of case studies showing that the existence of heterogeneous base states was not guaranteed.
Further, we could not determine, a priori, whether base states exist for a finite size heterogeneity.
To investigate this further, two parameter scans were performed.
The first varied the average production rate of the first chemical, and the length of the domain.
The second varied the form of the heterogeneity and the length of the domain.
Both parameter scans were tested with both the Schnakenberg and the Gierer-Meinhardt reactions.
For each set of parameters chosen, we measured how far the branch of solutions could be continued before reaching a fold bifurcation.
This measures how large the heterogeneity can be before the Turing base state ceases to exist.
The results of the parameter scans results reveal strong correlations with existing, fundamental theory from the dispersion relation.
Further research into a clear link between these theories is needed.

For small domain lengths, it becomes even more difficult to distinguish between patterned and non-patterned states.
This is because the wavelength of some patterns are often similar to the length scale of the heterogeneity.
The new definition allows for this distinction to be made, so systems with a small domain length can be analysed.
This new distinction allowed us to analyse how the critical domain length changes for heterogeneous RD systems.
We numerically determined the critical domain length for a range of heterogeneity sizes, and a range of average production rates.
This serves as a proof of concept of how the new definition could be applied to a new problem.
This was done for the Schnakenberg system and a Gierer-Meinhardt system.
We were able to find the critical domain length for a range of heterogeneity sizes and average production rates.
In some cases, however, the method we used to find the critical domain length failed.
It is possible that there are discontinuities in the critical domain length caused by a decoupling in the domain.
The method should be further developed to account for this, in an attempt to resolve the issues with the method used.



\bibliographystyle{siamplain}
\bibliography{references}
\end{document}